\documentclass{article} 

\usepackage{graphicx}
\usepackage{amsmath}
\usepackage{amsfonts}
\usepackage{latexsym}

\newtheorem{conjecture}{Conjecture}

\newtheorem{theorem}{Theorem}

\newtheorem{proposition}{Proposition}
\newtheorem{definition}{Definition}

\begin{document}

\title{\textbf{On certain positive integer sequences}}
\author{{\Large Giuseppe Melfi}\\
\textit{Universit\'e de Neuch\^atel}\\
\textit{Groupe de Statistique} \\
\textit{Espace de l'Europe 4}\\
\textit{CH-2002 Neuch\^atel}\\
\textit{Giuseppe.Melfi@unine.ch}
}
\date{}
\maketitle

{\abstract{A survey of recent results in elementary number
theory is presented in this paper. Special attention is given to
structure and asymptotic properties of certain families of positive
integers.
In particular, a conjecture on complete sequences
of Burr, Erd\H os, Graham and Wen-Ching Li is amended.}}

\section{Introduction}

The aim of this paper is to survey some topics in elementary number
theory. The contents are based on the talk given in Parma
at the occasion of the Second Italian Meeting of Number Theory
in november
2003.
All these topics are of arithmetical nature and, as is often
the case,
no special knowledge is required.
Still today improvements and new developments are possible, so this
makes this paper suitable for students who
wish to begin or to continue their studies in elementary number theory.

Each section contains a problem involving a suitable
sequence of family of sequences with a short literature, open problems,
and some new results. 

Section~\ref{prac} is consacrated to practical
numbers,
i.e., those numbers $m$ such that the set of all distinct positive
divisors sums contains all positive integers not exceeding~$m$.
Section~\ref{sumf} concerns sum-free sequences, i.e., increasing sequences of
positive integers such that each term is never a sum of distint
preceding terms. Section~\ref{sums} presents some results on density
of certain sets whose elements are sums of distinct powers of positive
integers. We provide a counterexample to a conjecture of Burr,
Erd\H os, Graham and Wen-Ching Li
\cite{Bu}, so showing the necessity of an amended version.
Section~\ref{simu} deals with certain families of positive
integer sequences whose digital expansion of elements is suitably
related to the digital expansion of their powers.
\section{Practical numbers}\label{prac}

A positive integer $m$ is a practical number if every positive integer $n<m$
is a sum of distinct positive divisors of $m$. 

Let $P(x)$ be the counting function of
practical numbers
and let $P_2(x)$ be the function that counts 
practical numbers $m\le x$ such that
$m+2$ is also a practical number. Stewart \cite{St} proved that 
a positive integer $\,m\ge 2,$
$\,m=q_1^{\alpha_1}q_2^{\alpha_2}\cdots q_k^{\alpha_k},$
with primes $q_1<q_2<\dots<q_k$ and integers $\alpha_i\ge 1,$
is practical\index{practical} if and only if $\,q_1=2\,$ and, for
$i=2,3,\dots,k,$
$$q_i\le\sigma\!\left(q_1^{\alpha_
1}q_2^{\alpha_2}\cdots q_{i-1}
^{\alpha_{i-1}}\right)+1,$$
where $\sigma(n)$ denotes the sum of the positive divisors of $n$.
A wide survey of results and conjectures on practical numbers is given 
by Margenstern~\cite{Ma2}.

Practical numbers appear to be a prime-like sequence.
Saias~\cite{Sai}, using suitable sieve methods introduced by Tenenbaum
 provided
a good estimate in terms of
a Chebishev-type theorem: 
\begin{theorem}
For suitable
constants $c_1$ and $c_2$,
$$c_1\frac x{\log x}<P(x)<c_2\frac x{\log x}.$$
\end{theorem}
The author \cite{Me1} solved the Golbach problem by proving that
every even positive integer is a sum of two practical numbers.
The proof used an auxiliary increasing sequence $m_n$ of practical numbers
such that for every $n$, $m_n+2$ is also a practical number and
$m_{n+1}/m_n$ bounded by an absolute constant, and a corollary of
Stewart's theorem, namely
if $m$ is a practical number and $n\le2m$, then $mn$ is a pratical
number too. Every pair of twin practical numbers yields a suitable
interval of even numbers expressible as a sum of two practical numbers
and intervals overlap.

A local property which does not appear in primes holds for practical numbers:
there exist infinitely many practical numbers 
$m$ such that  both $m-2$ 
and $m+2$ as one can check by taking $m=2\cdot3^{3^k\cdot70}$, for
$k\in\mathbb N$.

Further, infinitely many $5$-tuples of practical numbers of the form
$m-6$, $m-2$, $m$, $m+2$, $m+6$ exist under a suitable but
still unproved hypothesis  \cite{Me4}.

Twenty years ago Margenstern conjectured that 
for suitable $\lambda_1$ and $\lambda_2>0$
$P(x)\sim\lambda_1\frac x{\log x}$ and  
$P_2(x)\sim\lambda_2\frac x{(\log x)^2}.$
He empirically proposed $\lambda_1\simeq1.341$ and
$\lambda_2\simeq1.436$. Such conjectures appear far to be
proved. However
it should be interesting to prove that
$\lim_{x\rightarrow\infty}P(2x)/P(x)=2,$
a somewhat weaker result conjectured
by Erd\H os \cite{Er2}.

Concerning the counting function of twin practical numbers a recent
result of the author \cite{Me3} is the following
\begin{theorem}
Let $k>2+\log(3/2)$. For sufficiently large $x$, 
$$P_2(x)>\frac x{\exp\{k(\log x)^\frac12\}}.$$
\end{theorem}
In particular this implies that for every $\alpha<1$, \ $P_{2}(x)\gg
x^\alpha$. The proof uses the fact that if $m_1$ and $m_2$ are
two practical numbers, with $0.5<m_1/m_2<2$, and
with ${\rm gcd}\{m_1,m_2\}=2$, there exist $r$ and $s$, not exceeding
respectively
$2m_1$ and $2m_2$,
such that $m_1r$ and $m_2s$ are a pair of twin practical numbers.
It is possible to build $m_1$ and $m_2$ by pick primes of their
factorization
in suitable intervals in order to control their mutual size.
One has to count all possible pairs and divide by 
the maximal number of repetitions.

Many other open probems on practical numbers and related questions
have been raised by Erd\H os in \cite{Erdlast}.

\section{Sum-free sequences}\label{sumf}

An increasing sequence of positive integers $\{n_1,n_2,\dots\}$ 
is called a sum-free sequence if each term is never
a sum of distinct smaller terms.
This definition is due to Erd\H os \cite{Er3} who proved certain related 
results and raised several problems. In his paper he proved that
if $\{n_k\}$ is a sum-free\index{sum-free sequence} sequence
then it has zero asymptotic density. In other words, for every
$\varepsilon>0$, and for sufficiently large $k$, \ $n_k>k^{1+\varepsilon}$.
By the same argument he proved that for every
$\beta<(\sqrt5+1)/2$, for infinitely many $k$, \ $n_k>k^\beta$.

Until 1996, all known sum-free sequences had a gap, namely
$$\limsup_{k\rightarrow\infty}\frac{n_{k+1}}{n_k}>1.$$ 
In
\cite{dezuerdos},
Deshouillers, Erd\H os and the author
gave some examples of infinite 
sum-free sequences with no gap.
They also proved that for every positive $\delta$,
there exists a sum-free sequence 
$\{n_k\}$ such that $n_k\sim k^{3+\delta}.$

The best extremal result concerning sum-free sequences is due to 
\L uczak and Schoen \cite{LuSc}. They proved that
for every positive $\delta$,
there exists a sum-free sequence 
$\{n_k\}$ such that 
$n_k\sim k^{2+\delta}$ and 
that the exponent $2$ is the best possible. 

However an interesting problem remains open.
Erd\H os
proved that for any sum-free sequence $\{n_k\} $
one has
\[\sum_{j=1}^{\infty}\frac1{n_j}<103.\]
It is natural to define 
\[R=\sup_{\{n_k\} \mbox{ {\scriptsize sum-free}}}
\left\{\sum_{j=1}^{\infty}\frac1{n_j}
\right\}.\]
The best known lower bound is due to Abbott~\cite{Ab} who proved that 
$R>2.064$.
Levine and 
O'Sullivan~\cite{Le-O'S} proved the actual best known upper bound $R<4.$

\section{Sums of distinct powers}\label{sums}
A sequence $S=\{s_1,s_2,\dots\}$ of positive integers is a \textsl{complete}
sequence, if
$$\Sigma(S):=\left\{\sum_{i=1}^\infty\varepsilon_i s_i,\  \ \ \ \mbox{ for }
\varepsilon_i\in\{0,1\},\ \ \ \ 
\sum_{i=1}^\infty\varepsilon_i<\infty \right\}$$
contains all sufficiently large integers.
In this section we deal with sequences $S$ whose elements are powers of 
positive integers. 
Let $s\ge1$ and $A$ be a (finite or infinite) set of integers greater
than 1.
Let ${\rm{Pow }}(A;s)$ be the nondecreasing sequence of positive
integers 
of the form $a^k$ with $a\in A$ and $k\ge s$. 

Burr, Erd\H os, Graham and Wen-Ching Li~\cite{Bu} proved several results
providing sufficient conditions in order that Pow$(A;s)$ is complete. 
They also conjectured that
for any $s\ge1$, ${\rm{Pow }}(A;s)$ 
is complete
if and only if
\begin{itemize}
\item[(i)]$\sum_{a\in A}1/(a-1)\ge1$,
\item[(ii)] ${\rm{gcd }}\{a\in A\}=1$.
\end{itemize}
As Graham noted \cite{Gr}, $A$ must be intended as finite, even if from
the text of the conjecture this is not explicitely said. In fact
the following proposition disproves the conjecture in the case of 
infiniteness is allowed. 
We provide
counterexamples using suitable infinite sets $A$. However, we point out that  
for finite sets the problem is open.

\begin{proposition}
Let $\varepsilon\ge0$.
There exists a 
set $A$ of integers $\ge2$
such that:
\begin{itemize}
\item[(i)]$\sum_{a\in A}1/(a-1)<\varepsilon$,
\item[(ii)]  for every $s\ge1$, ${\rm{Pow }}(A;s)$ is complete.
\end{itemize}
\end{proposition}
\noindent 
\textbf{Proof.}
Let $p\ge3$ be a prime, and let $R_{p}:=\{n^2p,\ n\in\mathbb N\}$.
We have that  ${\rm{Pow }}(R_p;s)\supseteq\{n^{2s}p^s, \ n\in\mathbb N \}$.
Since every sufficiently large integer is a sum of distinct $2s$-th powers of
positive integers \cite{sp}, there exists $r_{2s}$ such that
$$\Sigma({\rm{Pow }}(R_p;s))\supseteq\{np^s, \ n>r_{2s}\}.$$

If $Q_p$ is such that 
$Q_p\cap R_p=\emptyset$, we get $\Sigma({\rm{Pow }}(R_p\cup Q_p;s))
=\Sigma({\rm{Pow }}(R_p;s))+\Sigma({\rm{Pow }}(Q_p;s))$. Since 
$\Sigma({\rm{Pow }}(R_p;s))$ contains all sufficiently large multiple of 
$p^s$, in order that ${\rm{Pow }}(R_p\cup Q_p;s)$ is complete,
it suffices to provide a set $Q_p$ such that 
$\Sigma({\rm{Pow }}(Q_p;s))$ contains at least one element for each congruence
class modulo $p^s$. 
Let $Q_p=\{p+1\}$. It is clear that ${\rm{Pow }}(Q_p;s)$ contains
infinitely many elements $\equiv1\bmod p^s$, so 
$\Sigma({\rm{Pow }}(Q_p;s))$ contains at least one element for each congruence
class modulo $p^s$.

By the above arguments this implies that for 
$A:=R_p\cup Q_p$, 
${\rm{Pow }}(A;s)$ is complete.
Note that elements of $A$ do not depend on $s$, and that for
sufficiently large $p$,
$$\sum_{a\in A}\frac1{a-1}=\frac1p+\sum_{n=1}^\infty\frac1{n^2p-1}<\varepsilon.$$
$\Box$
 
A related open question is the following one. Consider the sequence
$\{n_k\}$ of
positive integers that are a sum of distinct powers of 3 and of 4.
Erd\H os asked for a proof that $n_k\ll k$. The best known
result is $n_k\ll k^{1.0353}$ as shown in~\cite{Me2}.

More generally we propose the following conjecture.

\begin{conjecture}
Let $s\ge1$ and
let $A$ be a sequence of integers $\ge2$. If for every
$a_1,a_2\in A$, ${\rm gcd}\{a_1,a_2\}=1$ and 
$\sum_{a\in A}\frac1{\log a}>\log 2$
then $\Sigma({\rm{Pow }}(A;s))$ has positive lower asymptotic density.
\end{conjecture}

Note that if we replace the condition `for every
$a_1,a_2\in A$, ${\rm gcd}\{a_1,a_2\}=1$' by `${\rm gcd}\{a\in A\}=1$'
the statement is not true. For the set $A=\{3,9,81,104\}$, we have
${\rm gcd}\{a\in A\}=1$, and  $\Sigma({\rm{Pow }}(A;s))$ has
zero lower asymptotic density \cite{Me2}.

\section{Simultaneous binary expansions}\label{simu}

Let $B(n)$ be the sum of digits of the positive integer $n$
written on base~$2$. A natural question of some interest is to 
describe the sequence of positive integers $n$ such that
$B(n)=B(n^2)$. Other related questions can be easily raised. For example
it can be of interest to study the sequence the positive integers $n$ such that
$2B(n)=B(n^2)$, i.e., those $n$ such that $n$ and $n^2$ have the 
same `density' of ones in their binary expansion. 

\begin{definition}
Let $k\ge2$, $l\ge1$, $m\ge2$ be positive integers.
We say that a positive integer $n$ is a $(k,l,m)$-number if the 
sum of digits of $n^m$
in its expansion in base $k$ is $l$ times the sum of the digits of the
expansion in base $k$ of $n$.
\end{definition}  

The above sequences respectively represent the $(2,1,2)$-numbers and
the $(2,2,2)$-numbers.

The simplest case is $(k,l,m)=(2,1,2)$, which corresponds to the
positive integers $n$ for which the numbers of ones in their binary
expansion is equal to the number of ones in $n^2$.

The list of $(2,1,2)$-numbers as well as the list of $(2,2,2)$-numbers
shows several interesting facts. 
The distribution is not
regular.  A huge amount of questions, most of which of
elementary nature, can be raised. 

In spite of their elementary definition, as far as we know
these sequences do not appear in literature.
Several questions, concerning both the structure 
properties and asymptotic behaviour, can be raised. 
Is there a necessary and sufficient 
condition to assure that a number is of type (2,1,2)? 
of type (2,2,2)? What is the asymptotic
behaviour of the counting function of (2,1,2)-numbers? of (2,2,2)-numbers? 

The irregularity of distribution
does not suggest a clear answer to these questions.

Let $p_{(k,l,m)}(n)$ be the number of $(k,l,m)$-numbers 
which do not exceed~$n$. By elementary arguments one can prove the following.
\begin{theorem}
Let $p_{(2,1,2)}(n)$ be the counting function of the $(2,1,2)$-numbers. 
We have
$$p_{(2,1,2)}(n)\gg n^{0.025}.$$
\end{theorem}

The proof uses the fact the for every $n$ it is possible to construct
 a set of $n$ distinct (2,1,2)-numbers
not exceeding $An^{40}$. To do this, one uses the fact that for every
$n<2^\nu$, $B(n(2^\nu-1))=\nu$. The construction uses an arbitrary
number not exceeding $n$ and by adding in their binary expansion a
suitable finite sequence of zeros and ones, with a special attention
for the control of the function $B$ for simultaneously
the new number and its square, one obtains a $(2,1,2)$-number
not exceeding $An^{40}$. Further details can be found on~\cite{Me5}.

By an analogous procedure it is possible to prove the following theorem.

\begin{theorem}
Let $p_{(2,2,2)}(n)$ be the counting function of the $(2,2,2)$-numbers. 
We have
$$p_{(2,2,2)}(n)\gg n^{0.909}.$$
\end{theorem}

Concerning upper bounds for counting functions,
recently S\'andor~\cite{Sa} announced that
$$p_{(2,1,2)}(n)\ll n^{0.9183}.$$

Some conjectures on counting functions can be proposed.
Apart from small intervals centered in powers of 2, $B(n)$ and
$B(n^2)$ appear as a random sequence of zeros and ones. Using this
appearence and considering $B(n)$ and $B(n^2)$ as indipendent random
variables
where zeros and ones are equally probable, an eurisitc approach
suggests the following conjectures.  

\begin{conjecture}\label{conjecture}
$$p_{(2,1,2)}(n)= n^{\alpha+o(1)}$$
where $\alpha=\log{1.6875}/\log2\simeq0.7548875$.
\end{conjecture}

\begin{conjecture}\label{conjecture2}
For each $k$ one has:
$$p_{(2,k,k)}(n)= 
\frac{n}{(\log n)^{1/2}}G_k
+R(n),$$
where 
 $G_k=\sqrt{\frac{2\log2}{\pi(k^2+k)}}$ and
$R(n)=o(n/(\log n)^{1/2}).$
\end{conjecture}
A detailed discussion of the above conjecture can be found on~\cite{Me5}.
Computations show that the above conjectures describe quite well the
behaviour of counting functions for $n<10^8$.

\end{document}